\documentclass[dvips, 12pt,a4paper]{article}
\usepackage{latexsym, amsmath}
\usepackage{amssymb}
\begin{document}
\begin{center}
\textbf{Orbits of Real Algebraic Groups}\\[5mm]
H.Azad\\[4mm]
Department of Mathematics and Statistics\\
King Fahd University of Petroleum \& Minerals\\
Dhahran, Saudi Arabia\\
\texttt{hassanaz@kfupm.edu.sa}\\
\end{center}

\textbf{Abstract:} The existence of closed orbits of real algebraic groups on certain real algebraic spaces is established. As an application it is shown that if $G$ is a real reductive group with Iwasawa decomposition $G=KAN$, then all unipotent subgroups of $G$ are conjugate to a subgroup of $N$.\\[0.07in]
2000 MSC:  14L30, 20G20, 22E15, 17B81\\

\baselineskip=18pt

A real algebraic group $G$ is, for the purposes of this paper, a closed connected subgroup of $GL(n,R)$ such that the group $G^{\,C}$ of $GL(n,C)$whose Lie algebra is $Lie(G)+iLie(G)$ is an affine algebraic subgroup of $GL(n,C)$:here $Lie(G)$ is the Lie algebra of $G$.

The aim of this note is to prove the following result\\
{\bf Proposition:} $Let\, G\subset GL(n,R)$\emph{be a real algebraic group and H,L real algebraic subgroups of G.Let} $M = G^{\,C}/H^{\,C} and \,\,\, \xi_{0}=eH^{\,C}.$\emph{Then L has a closed orbit in}$\, G\,\xi_{0}.$

This result has an interesting corollary which is useful in several contexts.\\
{\bf Corollary:} $Let\ G$ \emph{be a reductive group with Iwasawa decomposition} $G=KAN$. Let
$u \subset Lie(G)$ \emph{be a subalgebra in which every element is nilpotent and U the subgroup of G whose Lie algebra is u. Then U is conjugate to a subgroup of N. In particular, all maximal unipotent subgroups are conjugate in G}

This result is also a special case of a result in Borel-Tits [2 p.126]. It is also proved in Mostow [5] and Vinberg [6 p.276, 11]. It seems not be as well known as it should be. For example, it is not given in standard references like Helgason  [3] or Knapp [4], although the 1 dimensional case is proved in Helgason [3 p.431] using the Jacobson-Morozov theorem. The author learnt about the results in [2] from Venkatramana (private communication) after a first draft of the note had been written.

A knowledge of solvable subgroups is of importance in theoretical physics, as explained in the papers of Patera, Winternitz and Zassenhaus [8,9], where the authors have determined all maximal solvable connected subgroups of the classical real groups. The classification of solvable subgroups is also of great practical use in the reduction theory of differential equations.

Our proof is based on the following ideas, whose details are given in the proof of Theorem 3 in [1]. For a subset $Z$ of a manifold $M$ and a point $p$ in $M$, the tangent space $T_{p}(Z)$ to $Z$ at $p$ is the subspace of $T_{p}(M)$ spanned by the vectors
$\displaystyle \gamma'(0),$ where $\gamma:I\rightarrow M$ is a differentiable curve whose trace lies in $Z$ with $\gamma(0)=p$. \\
Let $K$ be a compact group of transformations of $M$ whose fixed point set $M_{K}$ is nonempty. For a point $p$ in $M_{K}$, the set of all points $r$ in $M_{K}$ which can be joined to $p$ by curves $\{\gamma(t)\}_{t\in I}\subset M_{K}$  with dim $T_{\gamma(t)}(M_{K})$ = dim $T_{p}(M_{K})$ for all $ t\in I$ is both open and closed in $M_{K}$. Hence it is the component of $M_{K}$ passing through $p$ and this component is invariant under any connected group operating on $M_{K}$. It follows that if $G \supset H$ are Lie groups and $\sigma$ is an automorphism of finite order of $G$ with $\sigma (H)=H,$ then the $(G_{\sigma})^{o}$ - orbits in $(G/H)_{\sigma}$ are the connected components of $(G/H)_{\sigma}$.

For a closed connected subgroup $G$ of $GL(n,R)$ its complexification $G^{C}$ is the subgroup of $GLn,C)$ whose Lie algebra is $Lie(G)+iLie(G)$.\\[0.1in]
{\bf Proof of the Proposition:}\\[0.05in]
The space $M=G^{C}/H^{C}$ is a smooth algebraic variety. Let $\sigma$ denote complex conjugation on $G^{C}$ as well on $M$. The group $G^{C}$ is generated by the complex 1 parameter subgroups $exp(zX),z\in C,$ with $X\in Lie(G)$, so indeed $\sigma$ operates on $G^{C}$ and on $G^{C}/H^{C}$ by $\sigma (gH^{C})=\sigma (g) H^C$, $g\in G^{C}$.

Consider the $L^{C}$ orbit of $\xi_{0}=eH^{C}$ in $M=G^{C}\xi_{0}$. Because all groups involved are algebraic, the Zaraski closure of $L^{C} \xi_{0}$ contains $L^{C} \xi_{0}$ and lower dimensional orbits [7,p.19].  Let $Z$ be the union of all the lower dimensional orbits. We have \\
$\overline{L^{C} \xi_{0}}= L^{C} \xi_{0}\cup Z$ \\[0.1in]
Now $\overline{(L^{C} \xi_{0})_{\sigma}}\cap G\xi_{0}= ((L^{C}\xi_{0})_{\sigma}\cap G \xi_{0})\cup((Z)_{\sigma}\cap G\xi_{0})\qquad ({*})$\\[0.1in]
Case 1: Suppose, in equation $({*})$ on the right hand side, the second intersection is empty.\\
We then have \\[0.1in]
$\overline{(L^{C} \xi_{0})_{\sigma}}\cap G \xi_{0}= (L^{C} \xi_{0})_{\sigma}\cap G\xi_{0}$\\[0.1in]
From now on, the topological terms used are in the context of the classical topology. From the left hand side of the last equation, the connected component of $\xi_{0}$ in
$\overline{(L^{C} \xi_{0})_{\sigma}}\cap G \xi_{0}$ is a closed set in $G \xi_{0}$. From the right hand side, the connected component of $\xi_{0}$ in  ${(L^{C} \xi_{0})_{\sigma}}\cap G \xi_{0}$ is at least $L\xi_0$ and it cannot be more because the connected component of $\xi_0$ in $(L^{C}\xi_0)_\sigma$ is $L\xi_0$.\\
So, $L\xi_0$ is a closed subset of the orbit $G\xi_0$.\\[0.1in]
Case2: Suppose in equation $(*)$, the second intersection is nonempty. The orbit $G\xi_0$ is a union of $L$ orbits. So take a point $g\xi_0 = \eta _0$ in $Z$. As $g$ is in $G$ we have $G\xi_0=G\eta_0$ but now the Zariski closure of the orbit $L^C \eta_0$ has the open $L^C\eta_0$ orbit with lower dimension than $L^C\xi_0$. Therefore, repeating the argument as in Case 1 with $\xi_0$ replaced by $\eta_0$ we find, ultimately that for some $g$ in $G$, the orbit $Lg\xi_0$ is closed in $G\xi_0$. This completes the proof of the proposition.\\[0.2in]
{\bf Proof of the Corollary}\\[0.1in]
Let $G=KAN$ be an Iwasawa decomposition of $G$. The group $A$ is an algebraic torus of $G$ in the sense that $A^C$ is an algebraic torus of $G^C$. Let $H=AN$. So $H^C=A^C N^C$ and $(H^C)_\sigma =(A^C)_\sigma (N^C)_\sigma=AN$

As $A^C$ and $N^C$ are algebraic groups, so is $H^C$ and $M=G^C/H^C$ is a smooth algebraic variety. Let the group $U$ be as in the statement of the corollary. Since every element in the Lie algebra of $U$ is nilpotent, the group $U$ is a real algebraic unipotent subgroup of $G$.

Let  $\xi_0 = eH^C$. The orbit $G\xi_0$ is compact and, by the Proposition, the group $U$ has a closed orbit in $G\xi_0$. Let this orbit be $Ug\xi_0$. On the one hand this is a closed set in the compact set $G\xi_0$; on the other hand, being an orbit of a unipotent subgroup, it must be a cell. Therefore it must be a point. Hence $Ug\xi_0=g\xi_0$, which means that $g^{-1}Ug \subset H^C$.

But  $G \cap H^C =H$, so $g^{-1}Ug \subset H$. Since the only unipotent elements of $H$ are in $N$ we have $U\subset gNg^{-1}$. This completes the proof of the corollary.

To determine maximal connected solvable algebraic groups of $G$, one notices that the unipotent part is normalized by the semisimple elements in the solvable group. Specific detailed information is in the basic papers of Patera et al [8,9] and in Snobel-Winternitz [10].\\

\subsection*{Acknowledgements:} The author thanks  KFUPM for funding Research Project  IN101026.

\end{document}